\documentclass[a4paper,conference]{IEEEtran}

\usepackage{indentfirst}
\usepackage[T1]{fontenc}
\usepackage[dvips]{graphicx}
\usepackage[dvips]{epsfig}
\usepackage{amsmath}
\usepackage{amssymb}
\usepackage{verbatim}
\usepackage{tikz} 
\usepackage{bbm}
\usepackage{algorithmic}
\usepackage{algorithm}
\usepackage{array}
\usepackage{enumerate}

\newtheorem{definition}{Definition}

\newcommand{\qed}{\nobreak \ifvmode \relax \else
      \ifdim\lastskip<1.5em \hskip-\lastskip
      \hskip1.5em plus0em minus0.5em \fi \nobreak
      \vrule height0.5em width0.5em depth0.25em\fi}

\newcommand{\E}[1]{\mathbf{E}_{}\left[#1\right]}

\begin{document}

\title{Disaster Recovery in Wireless Networks:\\ A Homology-Based Algorithm}
\author{\IEEEauthorblockN{A. Vergne, I. Flint, L. Decreusefond and
    P. Martins}
  \IEEEauthorblockA{
    Institut Telecom, TELECOM ParisTech, LTCI\\
    Paris, France\\
    Email: avergne, flint, decreuse, martins@telecom-paristech.fr} }

\maketitle

\begin{abstract}
In this paper, we present an algorithm for the recovery of wireless
networks after a disaster. Considering a damaged wireless network,
presenting coverage holes or/and many disconnected components, we 
propose a disaster recovery algorithm which repairs the network. It
provides the list of locations where to put new nodes in order to
patch the coverage holes and mend the disconnected components. 
In order to do this we first consider the simplicial complex
representation of the network, then the algorithm adds supplementary
vertices in excessive number, and afterwards runs a reduction
algorithm in order to reach an optimal result. One of the novelty of
this work resides in the proposed method for the addition of
vertices. We use a determinantal point process: the Ginibre point
process which has inherent repulsion between vertices, and has
never been simulated before for wireless networks representation. 
We compare both the determinantal point process addition method with
other vertices addition methods, and the whole disaster recovery
algorithm to the greedy algorithm for the set cover problem.
\end{abstract}

\section{Introduction}
Wireless networks are present everywhere, must it be sensor networks
or cellular networks. Fields where wireless sensor
networks can be used range from battlefield surveillance to target
enumeration in agriculture and include environmental
monitoring. In most applications, the topology of the
network, such as its connectivity and its coverage, is a critical
factor. Cellular networks are used for radio communication, where
coverage is also a critical factor. Indeed the covered area is often
the main characteristic of a cellular network. However such networks
are not necessary built with redundancy and can be sensitive to
disasters.

In case of a disaster, a wireless network can be seriously damaged:
some of its nodes can be completely destroyed. Coverage holes can
appear resulting in no signal  for communication or no monitoring at
all of a whole area, connectivity can be lost between nodes.  
Paradoxically, reliable and efficient communication and/or monitoring
are especially needed in such situations. Therefore solutions for
damage recovery  for the coverage 
of wireless networks are much needed. Extensive research on the
coverage problem in wireless networks exists: we can cite location-based
\cite{Fang} and range-based \cite{zhang} methods. 
However, connectivity based schemes seem of greater interest since
they provide an exact mathematical description of coverage without any
geographical (location or distance) information. In \cite{Ghrist1},
the authors introduced the Vietoris-Rips complex, based on the
proximity graph of a wireless network, as a tool to compute its
topology. Coverage computation via simplicial homology comes down to
linear algebra computations. It is for instance used in
\cite{Ghrist2} 
as a tool for a network operator to evaluate the
quality of its network.  

In this paper, we present a homology based algorithm for disaster
recovery of wireless networks. We represent wireless networks with 
\u{C}ech simplicial complexes characterizing their coverage. Given a set
of vertices and their coverage radius, our algorithm first adds
supernumerary vertices in order to patch every existing coverage hole
and connect every components, then runs an improved version of the
reduction algorithm presented in \cite{infocom} in order to reach an
optimal result with a minimum number of added vertices. At the end, we
obtain the locations in which to put new nodes. For the
addition of new vertices, we first compared two usual methods
presenting low complexity: grid positioning and uniform positioning.
Then, we propose the use of a
determinantal point process: the Ginibre point process. This process
has the ability to create repulsion between vertices, and therefore
has the inherent ability to locate areas with low density of vertices:
namely coverage holes. Therefore using this process, we will optimally
patch the damaged wireless network. The use and simulation of determinantal point
processes in wireless networks is new, and it provides tremendous
results compared to classic methods.We finally compared our
whole distaster recovery algorithm performance to the classic recovery
algorithm performance: the greedy algorithm for the set cover problem.

This is the first algorithm that we know of that adds too many
vertices then remove them to reach an optimal result instead of adding
the exact needed number of vertices. This, first, allows flexibility
in the choice of the new vertices positions, which can be useful when
running the algorithm in a real life scenario. Indeed, in case of a
disaster, every location is not always available for installing new
nodes and preferring some areas or locations can be done with our
algorithm. The originality of our work lies also in the choice of the
vertices addition method we suggest. 
On top of flexibility, our algorithm provides a more reliable repaired
wireless network than other algorithms. Indeed, adding the exact
needed number of vertices can be optimal mathematically speaking but
it is very  sensitive to the adherence of the nodes positions chosen
by the algorithm. To compare our work to literature, we can
see that the disaster recovery problem can be viewed as a set cover
problem. It suffices to define the universe as the area to be covered
and the subsets as the balls of radii the coverage radii. Then the
question is to find the optimal set of subsets that cover the
universe, considering there are already balls centered on the existing 
vertices. A greedy algorithm can solve this problem as explained in
\cite{greedy}. We can see in \cite{epsilon-nets} that $\epsilon$-nets
also provide an algorithm for the set cover problem via a sampling of
the universe. We can also cite landmark-based routing, seen in
\cite{fps1}, which, using furthest point sampling,
provides a set of nodes for optimal routing that we can
interpret as a minimal set of vertices to cover an area.

The remainder of this paper is structured as follows: after a
section on related work we present the main idea of our disaster
recovery algorithm in Section \ref{sec_mai} using some definitions
from simplicial homology. Then in Section \ref{sec_add}, we compare
usual vertices addition methods. In Section \ref{sec_det}, we expose
the determinantal method for new vertices addition. Section \ref{sec_ra}
is devoted to the reduction algorithm description. Finally in Section
\ref{sec_per} we compare the performance of the whole disaster
recovery algorithm with the greedy algorithm for the set cover
problem. We conclude in Section \ref{sec_ccl}.

\section{Recovery in cellular networks}
The first step of recovery in cellular networks is the
detection of failures. The detection of the failure of a cell occurs when its
performance is considerably and abnormally reduced. In \cite{kaschub},
the authors distinguish three stages of cell outage: degraded,
crippled and catatonic. This last stage matches with the event of a
disaster when there is complete outage of the damaged cells.
After detection, compensation from other nodes can occur through relay assisted
handover for ongoing calls, adjustments of neighboring cell sizes via
power compensation or antenna tilt. In \cite{amirijoo}, the authors
not only propose a cell outage management description but also
describe compensation schemes. These steps of monitoring and
detection, then compensation of nodes failures are comprised under the
self-healing functions of future cellular networks.

In this work, we are interested in what happens when self-healing is
not sufficient. In case of serious disasters, the compensation from
remaining nodes and traffic rerouting might not be sufficient to provide service
everywhere. In this case, the cellular network needs a manual
intervention: the adding of new nodes to compensate the failures of
former nodes. However a traditional restoration with brick-and-mortar
base stations could take a long time, when efficient communication is
particularly needed. In these cases, a recovery trailer fleet of base
stations can be deployed by operators \cite{att}, it has been for
example used by AT\&T after 9/11 events. But a question remains:
where to place the trailers carrying the recovery base stations. An
ideal location would be adjacent to the failed node. However, these
locations are not always available because of the disaster, and the
recovery base stations may not have the same coverage radii than the
former ones. Therefore a new deployment for the recovery base stations has
to be decided, in which one of the main goal is complete coverage of the
damaged area. This becomes a mathematical set cover problem. It can
been solved by a greedy algorithm \cite{greedy}, $\epsilon$-nets
\cite{epsilon-nets}, or furthest point sampling \cite{fps1}. But
these mathematical solutions provide an optimal mathematical result
that do not consider any flexibility at all in the choosing of the new
nodes positions, and that can be really sensitive to imprecisions in
the nodes positions.

\section{Main idea}
\label{sec_mai}
When representing a wireless network, one's first idea will be
a geometric graph, where sensors are represented by vertices, and an
edge is drawn whenever two sensors can communicate with each
other. However, the graph representation has some limitations; first
of all there is no notion of coverage. Graphs can be generalized to
more generic combinatorial objects known as simplicial
complexes. While graphs model binary relations, simplicial complexes
represent higher order relations. A simplicial complex is a
combinatorial object made up of vertices, edges, triangles,
tetrahedra, and their $n$-dimensional counterparts. Given a set of
vertices $V$ and an integer $k$, a $k$-simplex is an unordered subset
of $k+1$ vertices $[v_0,v_1\dots, v_k]$ where $v_i\in V$ and
$v_i\not=v_j$ for all $i\not=j$. Thus, a $0$-simplex is a vertex, a 
$1$-simplex an edge, a $2$-simplex a triangle, a $3$-simplex a
tetrahedron, etc. 

Any subset of vertices included in the set of the $k+1$ vertices of a
$k$-simplex is a face of this $k$-simplex. Thus, a $k$-simplex has
exactly $k+1$ $(k-1)$-faces, which are $(k-1)$-simplices. For example,
a tetrahedron has four $3$-faces which are triangles. An abstract simplicial
complex is a collection of simplices which is closed with respect to
the inclusion of faces, i.e. all faces of a simplex are in the set of
simplices. For details about
algebraic topology, we refer to \cite{hatcher}.

We consider as inputs the set of existing vertices: the nodes of a
damaged wireless network, and their coverage radii. We also need a
list of boundary nodes, these nodes can be fictional, but we need to
know the whole area that is to be covered.
We restrict ourselves to wireless networks with a fixed
communication radius $r$, but it is possible to build the \u{C}ech
complex of a wireless network with different coverage radii
using the intersection of different size coverage balls.
The construction of the \u{C}ech abstract simplicial complex for a
fixed radius $r$ is given:
\begin{definition}[\u{C}ech complex]
  Given $(X,d)$ a metric space, $\omega$ a finite set of $N$ points in
  $X$, and $r$ a real positive number. The
  \u{C}ech complex of  $\omega$, denoted
  $\mathcal{C}_r(\omega)$, is the abstract simplicial complex
  whose $k$-simplices correspond to $(k+1)$-tuples of vertices in
  $\omega$ for which the intersection of the $k+1$ balls of radius
  $r$ centered at the $k+1$ vertices is non-empty.
\end{definition}
The \u{C}ech complex characterizes the coverage of the wireless network. 
The $k$-th Betti numbers of an abstract simplicial complex $X$ are
defined as the number of $k$-th dimensional holes in $X$ and are
computed via linear algebra computations. For example, $ \beta_0$
counts the number of $0$-dimensional holes, that is the number of
connected components. And $\beta_1$ counts the number of holes in the
plane. Therefore the Betti number $\beta_1$ of the \u{C}ech complex
counts the number of coverage holes of the wireless network it
represents. 

The algorithm begins by adding new vertices in addition to the set of
existing vertices presenting coverage holes. We suggest here the use
of two common methods, and the new determinantal addition method. As
we can see in Section \ref{sec_add}, it is possible to consider
deterministic or random based vertices addition methods:
flexibility is one of the greatest advantage of our algorithm. In
particular, it is possible to consider a method with pre-defined
positions for some of the vertices in real-life scenarios.
For any non-deterministic method, we choose that the number of added
vertices, that we denote by $N_a$, is determined as follows. First, it is set to be the
minimum number of vertices needed to cover the whole area minus
the number of existing vertices: $N_a:=\lceil \frac{a^2}{\pi r^2} \rceil - N_i$, $a^2$ being the
  area to cover. This way, we take into account the
number of existing vertices $N_i$. Then the Betti
numbers $\beta_0$ and $\beta_1$ are computed via linear algebra thanks
to the simplicial complex representation. If there is still more than
one connected component, and coverage holes, then the number of added
vertices $N_a$ is incremented with a random variable $u$ following an
exponential growth: $N_a=N_a+u$, and $u=2u$.

The next step of our approach is to run the coverage reduction
algorithm from \cite{infocom} which
maintains the topology of the wireless network: the
algorithm removes vertices from the simplicial complex
without modifying its Betti numbers. At this step, we remove some of
the supernumerary vertices we just added in order to achieve an
optimal result with a minimum number of added vertices. 
We give in Algorithm \ref{alg_dra} the outline of the algorithm. 
\begin{algorithm}[h]
  \caption{Disaster recovery algorithm}
  \label{alg_dra}
  \begin{algorithmic}[H]
    \REQUIRE Set of vertices $\omega_i$, radius $r$,
    boundary vertices $L_b$
\STATE Computation of the \u{C}ech complex $X=\mathcal{C}_{r}(\omega_i)$\;
\STATE  $N_a=\lceil \frac{a^2}{\pi r^2} \rceil - N_i$\;
\STATE Addition of $N_a$ vertices to $X$ following chosen method\;
\STATE Computation of $\beta_0(X)$ and $\beta_1(X)$\; 
\STATE $u=1$\;
\WHILE {$\beta_0\neq 1$ or $\beta_1 \neq 0$}
\STATE $N_a=N_a+u$\;
\STATE $u=2*u$\;
\STATE Addition of $N_a$ vertices to $X$ following chosen method\;
\STATE Computation of $\beta_0(X)$ and $\beta_1(X)$\; 
\ENDWHILE
\STATE Coverage reduction algorithm on $X$\;
\RETURN List $L_a$ of kept added vertices.
  \end{algorithmic}
\end{algorithm}

\section{Vertices addition methods}
\label{sec_add}
In this section, we propose two vertices addition methods. The aim of
this part of the algorithm is to add enough vertices to patch the coverage of the simplicial
complex, but the less vertices the better since the results will be
closer to the optimal solution. We consider grid and uniform positioning which
 require minimum simulation capacities and are well known in wireless
 networks management. 

The first method we suggest is deterministic: the number of added
vertices and their positions are set and are independent from the
initial configuration. It thus ensures perfect coverage,
the new vertices are positioned along a square lattice grid of
parameter $\sqrt{2}r$. 
In the second method we propose, the number of added vertices $N_a$ is computed accordingly to
the method presented in Section \ref{sec_mai}, taking into account the
number of existing vertices $N_i$. Then the $N_a$ vertices are sampled
following a uniform law on the entire domain. 

 We can compare the vertices addition methods presented here along two
variables: their complexity and their efficiency. 
First, we compare the complexities of the two methods. They
both are of complexity $O(N_a)$: computations of $N_a$ positions.
For the uniform method we have to add the complexity of computing the
coverage via the Betti numbers,
which is of the order of the number of triangles times the number of
edges that is $O((N_a+N_i)^5(\frac{r}{a})^6)$ for a square of side $a$
according to \cite{chimean}.    

To compare the methods
efficiency we count the number of vertices each method adds on average
to reach complete coverage. The grid
method being determinist, the number of added vertices is constant:
$N_a=(\lfloor \frac{a}{\sqrt{2}r} \rfloor +1)^2$ for a \u{C}ech
complex or $N_a=(\lfloor \frac{a}{2r} \rfloor +1)^2$ for a
Vietoris-Rips complex which is an approximation of the \u{C}ech complex
easier to simulate. We
can see in Table\ref{tableadded} the mean number of added vertices on
$1000$ simulations for each method in different scenarios on a square
of side $a=1$ with coverage radius $r=0.25$,  and a Vietoris-Rips
complex. Scenarios are defined by the mean percentage of area
covered before running the recovery algorithm: if there are many or
few existing vertices, and thus few or many vertices to add. We need to note that
number of added vertices is computed following our incrementation method
presented in Section \ref{sec_mai} and these results only concern the
vertices addition methods before the reduction algorithm runs.
The grid method is mathematically optimal for
the number of added vertices to cover the whole area, however it is
not optimal in a real life scenario where positions can not be defined
with such precision, and any imprecision leads to a coverage hole.  
This method fares even or better both in complexity and in number of
added vertices.

\section{Determinantal addition method}
\label{sec_det}
The most common point process in
wireless network representation is the Poisson point process. However
in this process, conditionally to the number of vertices,
their positions are independent from each other (as in the uniform
positioning method presented in Section \ref{sec_add}). This
independence creates some aggregations of vertices, that is not
convenient for our application. That is why we introduce the use of
determinantal point processes, in which the vertices positions are not
independent anymore. General point processes can be characterized by
their so-called Papangelou intensity. Informally speaking, for $x$ a
location, and $\omega$ a realization of a given point process, that is
a set of vertices, $c(x,\omega)$ is the probability to have a vertex
in an infinitesimal region around $x$ knowing the set of vertices
$\omega$. For Poisson process, $c(x, \omega)=1$ for any $x$ and any
$\omega$. A point process is said to be repulsive (resp. attractive)
whenever  $c(x,\omega) \ge c(x, \zeta)$ (resp. $c(x, \omega) \le c(x,
\zeta)$) as soon as $\omega \subset \zeta$. For repulsive point
process, that means that the greater the set of vertices, the smaller the probability
to have an other vertex.

Among repulsive point processes, we are in particular interested in
determinantal processes:
\begin{definition}[Determinantal point process]
Given $X$ a Polish space equipped with the Radon measure $\mu$, and $K$
a measurable complex function on $X^2$, we say that $N$ is a
determinantal point process on $X$ with kernel $K$ if it is a point
process on $X$ with correlation functions $\rho_n(x_1,\dots,x_n)=\det
(K(x_i,x_j)_{1\leq i,j \leq n})$ for every $n\geq 1$ and $x_1,\dots
,x_n \in X$.
\end{definition}
We can see that when two vertices $x_i$ and $x_j$ tends to be close to
each other for $i\neq j$, the
determinant tends to zero, and so does the correlation function. That
means that the vertices of $N$ repel each other.
There exist as many determinantal point processes as
functions $K$. We are interested in the following:
\begin{definition}[Ginibre point process]
The Ginibre point process is the determinantal point process with kernel
$K(x,y)=\sum_{k=1}^{\infty}B_k\phi_k(x)\overline {\phi_k(y)}$, where
$B_k, k=1,2,\dots$, are $k$ independent Bernoulli variables and
$\phi_k(x)=\frac{1}{\sqrt{\pi k!}}e^{\frac{-|x|^2}{2}}x^k$ for $x \in
\mathbb{C}$ and $k \in \mathbb{N}$. 
\end{definition}

The Ginibre point process is invariant with respect to
translations and rotations, making it relatively easy to simulate on a
compact set. Moreover, the repulsion induced by a Ginibre point
process is of electrostatic type. The principle behind the repulsion
lies in the probability density used to draw vertices positions. The
probability to draw a vertex at the exact same position of an already
drawn vertex is zero. Then, the probability increases with increasing
distance from every existing vertices. Therefore the probability to
draw a vertex is greater in areas the furthest away from every
existing vertices, that is to say in coverage holes. Therefore, added
vertices are almost automatically located in coverage holes thus
reducing the number of superfluous vertices.

Using determinantal point processes allows us to not only take into
account the number of existing vertices, via the computation of
$N_a$, but also take into account the existing vertices positions,
then every new vertex position as it is added. It suffices to consider
the $N_i$ existing vertices as the $N_i$ first vertices sampled in the
process, then each vertex is taken into account as it is drawn.
The Ginibre process is usually defined on the whole plane thus we
needed to construct a process with the same repulsive characteristics
but which could be restricted to a compact set. We also needed to
be able to set the number of vertices to draw. Due to space
limitations, we will not delve into these technicalities but they are
developed in \cite{ian}.

We can compare the determinantal vertices addition method to the
methods presented in Section \ref{sec_add}.
As for the complexity, since the determinantal method takes into
account the position of both existing vertices and randomly added
vertices, it is the more complex. First taking into account the
existing vertices positions is of complexity $O(N_i^2)$, then the
position drawing with the rejection sampling is of complexity
$O(N_a(N_a+N_i))$ at most. Thus we have a final complexity of
$O(N_i^2+N_a^2+N_aN_i))$. To which we add the Betti numbers
computation complexity: $O((N_a+N_i)^5(\frac{r}{a})^6)$.
We also give the comparison between the mean number
of added vertices for the three methods, simulation parameters
being the same as in Section \ref{sec_add}. The determinantal method
is the best method among all for the almost covered scenario. 
\begin{table}[h]
\centering
\begin{tabular}{|ccccc|}
\hline
\% of area initially  covered& 20\% & 40\%& 60\%& 80\%\\
\hline
Grid method&$9.00$&$9.00$&$9.00$&$9.00$\\
\hline
Uniform method&$32.51$&$29.34$&$24.64$&$15.63$\\
\hline
Determinantal method&$16.00$&$14.62$&$12.36$&$7.79$\\
\hline
\end{tabular}
\caption{Mean number of added vertices}
\label{tableadded}
 \end{table}
\vspace{-0.8cm}

\section{Reduction algorithm}
\label{sec_ra}
\begin{figure}[H]
      \scalebox{0.29}{\includegraphics{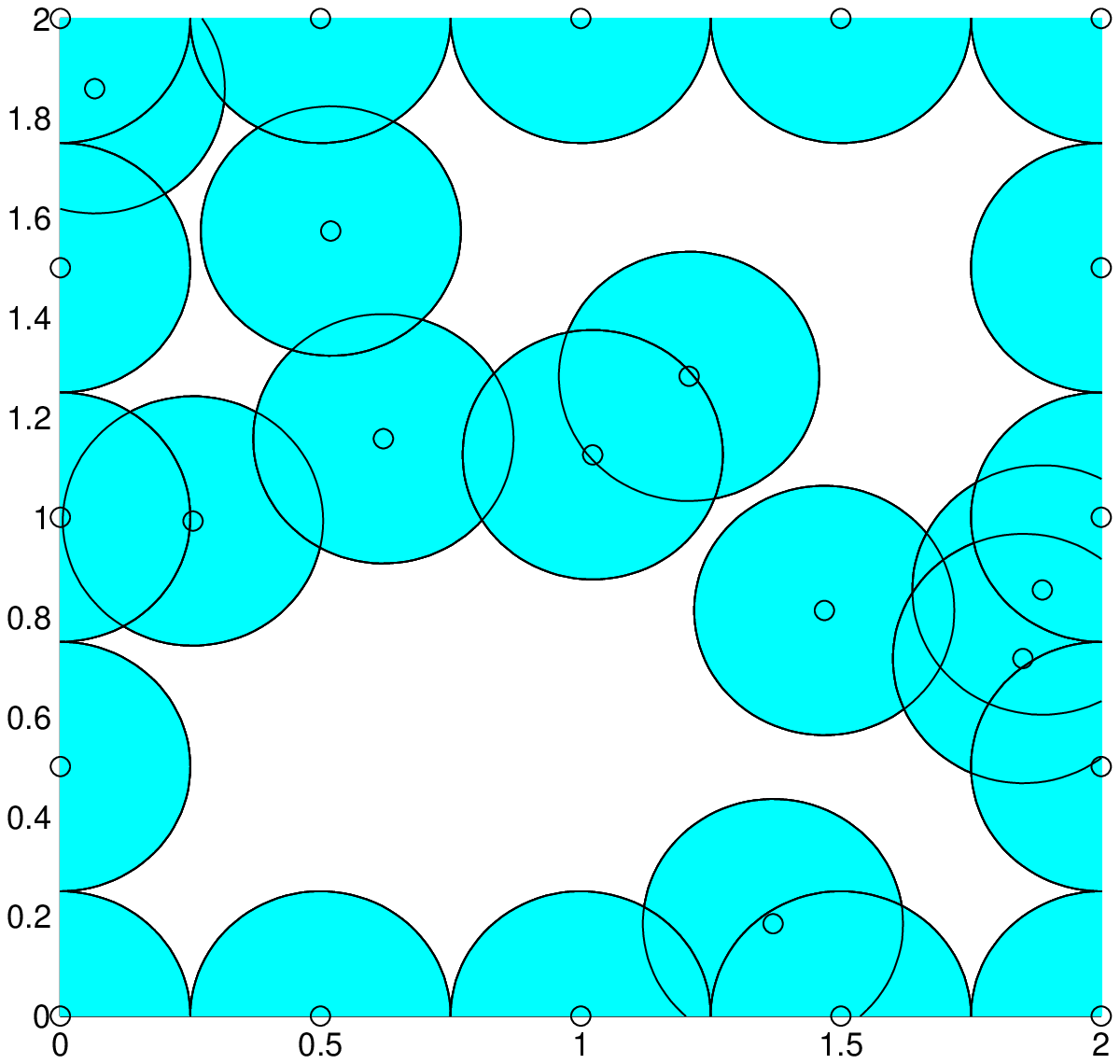}} 
\hfill
      \scalebox{0.29}{\includegraphics{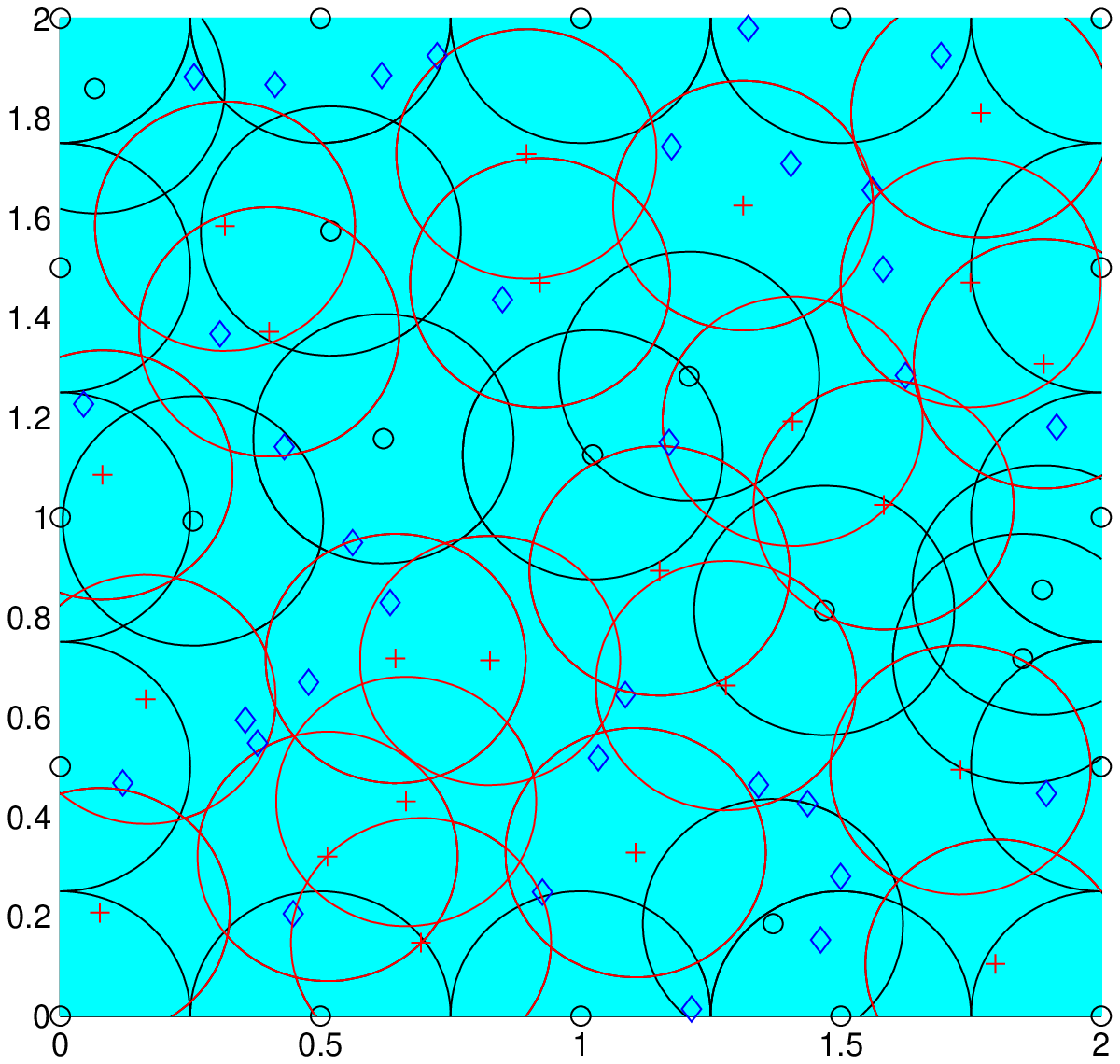}}
\caption{Execution of the homology algorithm}
    \label{fig_tot}
  \end{figure}

In this section, we recall the steps of the reduction algorithm
for simplicial complexes presented in \cite{infocom}. The algorithm takes
as input an abstract simplicial complex: here it is the \u{C}ech complex
of the wireless network plus the added vertices, and a list of
boundary vertices. At this step we have ensured that we have one 
connected component $\beta_0=1$, and no coverage hole $\beta_1=0$.
The first step is to characterize the superfluousness of
$2$-simplices for the coverage with a degree defined to be the size of
the largest simplex a $2$-simplex is the face of. Then to transmit the
superfluousness of its $2$-simplices to a vertex, an index of a vertex
is defined to be the minimum of the degrees of the $2$-simplices it is
a face of. The indices give an order for an optimal removal of vertices: the
greater the index of a vertex, the more likely it is 
superfluous for the coverage of its \u{C}ech simplicial complex. In
our disaster recovery case, we do not want to remove the remaining
vertices of our damaged network. So these remaining vertices are given
a negative index to flag them as unremovable,
and only the newly added vertices are considered for removing. So,
the vertices with the greatest index are candidates for removal: one
is chosen randomly. If its removal does not change the homology, then
it is effectively removed, otherwise it is flagged as unremovable with
a negative index. The algorithm goes on until every remaining vertex
is unremovable, thus achieving optimal result. For more information on
the reduction algorithm we refer to \cite{infocom}.
We can see in the Fig. \ref{fig_tot} an execution of the homology based
disaster recovery algorithm on a damaged network.
Existing vertices are black circles, kept added vertices are red
plusses, and removed vertices are blue
diamonds. 

\section{Performance comparison}
\label{sec_per}
We now compare the performance results of the whole disaster recovery
algorithm to the most known coverage recovery algorithm: the greedy
algorithm for the set cover problem. First, we compare their
complexities. The greedy algorithm method lays 
a lattice grid of parameter $\sqrt{2}r$ of
potential new vertices. Then the first added vertex is the furthest
from all existing vertices. The algorithm goes on adding the furthest
potential vertex of the grid from all vertices (existing+added). 
It stops when the furthest vertex is in the coverage ball of an
existing or added vertex. Then for the $(i+1)$-th vertex
addition, the greedy algorithm computes the distances from the
$N_i+i$ existing vertices to the $(\lfloor \frac{a}{\sqrt{2}r}
\rfloor +1)^2 - i$ potential vertices. Therefore its complexity
is in $O((N_i+N_a)(\lfloor \frac{a}{\sqrt{2}r} \rfloor +1)^2))$.  
For our  algorithm, we consider first the complexity
of building the abstract simplicial complex of the network which
is in $O((N_i+N_a)^C)$, where $C$ is the clique number. This
complexity seems really high since $C$ can only be upper bounded by
$N_i+N_a$ in the general case but it is the only way to compute the
coverage when vertices position are not defined along a given pattern. Then the
complexity of the coverage reduction algorithm is in
$O((1+(\frac{r}{a})^2)^{N_i+N_a})$ (see \cite{infocom}).  So
the greedy algorithm is less complex than ours in the general
case, however when $r$ is small before $a$, the trend is reversed. 

Then, we compare the mean number of added vertices in the final state.
The final number of added
vertices is  the number of added vertices for the greedy algorithm,
and the number of kept added vertices after the reduction
for the homology algorithm. It is important to note
that our algorithm with the grid method gives the exact same result
as the greedy algorithm, number of added vertices and their positions
being exactly the same. Simulations are done in the same conditions as
in Section \ref{sec_add}.
\begin{table}[h]
\centering
\begin{tabular}{|ccccc|}
\hline
\% of area initially  covered& 20\% & 40\%& 60\%& 80\%\\
\hline
Greedy algorithm&$3.69$&$3.30$&$2.84$&$1.83$\\
\hline
Homology algorithm&$4.42$&$3.87$&$2.97$&$1.78$\\
\hline
\end{tabular}
\caption{Mean final number of added vertices $\E{N_f}$}
\label{tablefinal}
 \end{table}
\vspace{-0.8cm}
The mean numbers of vertices added in the final state both with our recovery
algorithm and the greedy algorithm are presented in Table
\ref{tablefinal}. They are roughly the same, they both
tend to the minimum number of vertices required to cover the uncovered
area depending on the initial configuration. Nonetheless, we can see
that our algorithm performs a little bit worse than the greedy
algorithm in the less covered area scenarios because the vertices are
not optimally positioned and it can be seen when just a small
percentage of area is covered, and whole parts of the grid from the
greedy algorithm are used, instead of isolated vertices. In
compensation, our homology algorithm performs better in more covered
scenarios. The greedy algorithm is not flexible
at all: its success depends highly on the precision of the chosen
positions: a coverage hole appears as soon as a vertex is slightly
moved. Therefore our algorithm seems more fitted to the disaster
recovery case when a recovery network is deployed in emergency
both indoor, via Femtocells,  and outdoor, via a trailer fleet, where
GPS locations are not always available, and sticking to
positioning not always exact.

 \section{Conclusion}
\label{sec_ccl}
 In this paper, we adopt the simplicial homology representation for
 wireless networks which characterizes both the connectivity and
 the coverage of a given network. Based on that representation,
 we write an algorithm which patches coverage holes of damaged
 wireless networks by giving the positions in which to put new nodes.
Our recovery algorithm first adds enough new nodes to cover the whole
domain, then runs a reduction algorithm on the newly added nodes to
reach an optimal result. The originality of the algorithm lies in the
fact that we do not only add the needed nodes, thus providing a
mathematically optimal but not reliable result, but adds too many
nodes before removing the superfluous ones thus providing a stronger coverage that is less sensitive to the
imprecisions of following approximatively GPS locations. Moreover, the
vertices addition methods can be adapted to any particular
situation. We compare classic positioning methods to the new
determinantal method that is more efficient at 
positioning new vertices where they are needed.

\bibliographystyle{abbrv} \bibliography{article}

\end{document}